\RequirePackage{ifpdf}

\newif\iffinalversion
\finalversiontrue

\newif\ifnormalpage
\normalpagetrue

\ifnormalpage
  \ifpdf
    \iffinalversion
      \documentclass[pdftex,a4paper,notitlepage]{article}
    \else
      \documentclass[pdftex,a4paper,notitlepage]{article}
    \fi
  \else
    \iffinalversion
      \documentclass[hypertex,a4paper,notitlepage]{article}
    \else 
      \documentclass[hypertex,a4paper,notitlepage]{article}
    \fi
  \fi
\else
  \ifpdf
    \documentclass[pdftex,a3paper,twocolumn,notitlepage]{article}
    \usepackage[landscape,foot=1cm,head=1cm,rmargin=2cm,lmargin=2cm,twoside=false,mag=1000]{geometry}
  \else
    \documentclass[hypertex,a4paper,notitlepage]{article}
    \usepackage[landscape,foot=1cm,head=1cm,rmargin=2cm,lmargin=2cm,twoside=false,mag=1000]{geometry}
  \fi
\fi
\iffinalversion
  \usepackage[usenoborder,usenoshowkeys,usehyperref,usenocolorlink]{load_packages}
\else
  \usepackage[useshowkeys,usehyperref]{load_packages}
\fi
\iffinalversion\else
\fi
\usepackage{url}
\usepackage{mystyle}
\ifpdf 
\fi
\iffinalversion\else
\fi
\newcommand{\Dedekind}[1]{\ensuremath{{#1}^{\mathcal D}}}
\newcommand{\DNo}{\ensuremath{\Dedekind{\No}}}

\newcommand{\DKfield}{\ensuremath{\Dedekind \Kfield}}
\newcommand{\BNo}{\ensuremath{\overline{\No}}}
\newcommand{\clos}[2]{\ensuremath{{#1}^{#2}}}
\newcommand{\closA}[1]{\ensuremath{\clos {#1}{\Afam}}}
\newcommand{\Afamf}{\ensuremath{\Afam\cup\{\ff\}}}
\newcommand{\closAf}[1]{\ensuremath{\clos {#1}{\Afam\cup\{\ff\}}}}

\newcommand{\function}[1]{\ensuremath{\mathpzc{#1}}}
\newcommand{\ff}{\ensuremath{\function{f}}}

\newcommand{\gf}{\ensuremath{\function{g}}}

\newcommand{\pf}{\ensuremath{\function{p}}}
\newcommand{\qf}{\ensuremath{\function{q}}}
\newcommand{\xv}{\vettore{x}}
\newcommand{\yv}{\vettore{y}}
\newcommand{\av}{\vettore{a}}

\newcommand{\dv}{\vettore{d}}
\newcommand{\ev}{\vettore{e}}

\newcommand{\lt}{\ensuremath{\ell}}

\newcommand{\simpleq}{\ensuremath{\preceq}}
\newcommand{\simple}{\ensuremath{\prec}}
\newcommand{\nsimpleq}{\ensuremath{\npreceq}}
\providecommand{\nsimeq}{\ensuremath{\mathrel{%
      \makebox[.9em]{$\simeq\kern-.9em/$}%
    }}}
\newcommand{\stless}{\ensuremath{\mathbin{\lnapprox}}}

\newcommand{\Succ}{\ensuremath{\mathcal{S}}}
\newcommand{\isom}{\ensuremath{\simeq}}

\newcommand{\cut}[2]{\ensuremath{\langle\,#1 \mid #2 \,\rangle}}

\newcommand{\Cut}[2]{\ensuremath{\bigl\langle\,#1 \bigm| #2 \,\bigr\rangle}}
\newcommand{\bigcut}[2]{\ensuremath{\Bigl\langle\,#1 \Bigm| #2 \,\Bigr\rangle}}

\newcommand{\Splitcut}[2]{\ensuremath{\bigl\langle\,#1 \bigm|\\%
    & \qquad #2 \,\bigr\rangle}}
\newcommand{\cuts}[2]{\ensuremath{(\,#1\mid#2\,)}}

\newcommand{\bigcuts}[2]{\ensuremath{\Bigl(\,#1\Bigm|#2\,\Bigr)}}

\newcommand{\concat}{\ensuremath{\mathbin{:}}}

\newcommand{\opt}[1]{\ensuremath{#1^{\,\mathbf o}}}
\newcommand{\Opt}[1]{\opt{\big(#1\big)}}

\newcommand{\xo}{\ensuremath{\opt x}}
\newcommand{\yo}{\ensuremath{\opt y}}
\newcommand{\ao}{\ensuremath{\opt a}}

\newcommand{\Do}{\ensuremath{\opt d}}

\newcommand{\cLp}{\ensuremath{{c^L}'}}
\newcommand{\cRp}{\ensuremath{{c^R}'}}

\newcommand{\avo}{\ensuremath{\opt{\av}}}
\newcommand{\ffo}{\ensuremath{\opt {\ff}}}

\newcommand{\pfo}{\ensuremath{\opt {\pf}}}

\newcommand{\ffL}{\ensuremath{{\ff}^L}}
\newcommand{\ffR}{\ensuremath{{\ff}^R}}

\author{Antongiulio Fornasiero}
\title{Recursive definitions on surreal numbers}

\begin{document}

\maketitle

\begin{abstract}
Let $\No$ be Conway's class of surreal numbers.
I will make explicit the notion of a function $\ff$ on $\No$ recursively defined over some family of functions.
Under some `tameness' and uniformity conditions, $\ff$ must satisfy some interesting properties; in particular, the supremum of the class
$
\set{x\in\No: \ff(x)\geq 0}
$
is actually an element of $\No$.
As an application, I will prove that concatenation function $x\concat y$ cannot be defined recursively in a uniform way over polynomial functions.
%
\end{abstract}
%

\tableofcontents

\section{Introduction}
The class of surreal numbers $\No$ was introduced by Conway in~\cite{conway}.
I will present some results on regarding the properties of functions defined recursively on $\No$.
As a particualar case, I will give a different proof of the fact that $\On$ is a real closed field.
%

I will assume on the part of the reader some familiarity with the theory of surreal numbers, as exposed in~\cite{conway,gonshor,alling};
however, I will repeat some of its fundamental properties along the way.
\subsection{Basic definitions and properties}
I will recall some of the basic properties of the class of surreal numbers.
The definitions and theorems of this section can be found in~\cite{conway} and~\cite{gonshor}.

I will work in the set theory \NBG of von Neumann, Bernays and G\"odel with global choice.
A well formed formula of \NBG is a formula with set and class variables, without quantifications over classes.

I remind that $\On$ is the class of all ordinals.
\begin{definizione}[Surreal numbers]\label{def:surreal}
Following Gonshor, I define a surreal number $x$ as a function with domain an ordinal $\alpha$ and codomain the set $\set{+,-}$.
The ordinal $\alpha$ is called the \intro{length} of the $x$, in symbol $\lt(x)$.
The collection of all surreal numbers is the proper class $\No$.
\end{definizione}
\begin{definizione}[Linear order]
On $\No$ there is a linear order, defined according to the rule
$x\leq y$~($x$ is less or equal to $y$) iff
\[
\begin{array}{r@{\,}c@{\,}lcr@{\,}c@{\,}ll}
x         &= &y &\text{or}\\
x(\gamma) &\multicolumn{2}{c}{\text{\ undefined}} &\text{and} &y(\gamma) &= &+ &\text{or}\\
x(\gamma) &= &- &\text{and} &y(\gamma) &= &+ &\text{or} \\
x(\gamma) &= &- &\text{and} &y(\gamma) & \multicolumn{3}{c}{\text{undefined},} \\[.5ex]
\multicolumn{4}{l}{\text{where}}\\[.5ex]
\gamma &:= &\multicolumn{6}{@{\,}l}{\min\set{\beta\in\On: x(\beta)\neq y(\beta)}.}
\end{array}\]
\end{definizione}
\begin{definizione}[Simpler]
There is also a partial order $x\simple y$ (x is strictly \intro{simpler} than $y$, or $x$ is a \intro{canonical option} of $y$),
iff $x$ is the restriction of $y$ to an ordinal strictly smaller than $\lt(y)$.
\end{definizione}
\begin{definizione}[Convex]
A subclass $A\subseteq \No$ is \intro{convex} iff
\[
\forall x,y\in A\ \forall z\in\No\quad x<z<y \to z\in\No.
\]
Given $L,R\subseteq \No$ and $z\in\No$, $L<z$ means that
\[
\forall x\in L\quad x<z;
\]
similar definition for $x<R$.\\
$L<R$ means that
\[
\forall x\in L\ \forall y\in R\quad x<y.
\]
\end{definizione}
The fundamental properties of $\simple$ and of $<$ are given by the following theorem.
\begin{thm}
\begin{itemize}
\item $(\No,\simpleq)$ is a well founded partial order.
\item $(\No,\leq)$ is a dense linear order.
\item $\forall a\in\No$, $\Succ(a) := \set{x\in\No: a\simpleq x}$ is a convex subclass of $\No$.
\item If $A\subseteq \No$ is convex and non-empty, then there is a unique simplest element $a\in A$, \ie 
\[
\forall x\in A\ a\simpleq x.
\]
\item If $L<R$ are sub\emph{set}s of $\No$, then the cut
\[
\cuts{L}{R}:= \set{x\in\No: L<x<R}
\]
is convex and non-empty.
Its simplest element is called $\cut{L}{R}$.
\item $\cut{L}{R} = \cut{L'}{R'}$ iff for every $x^L \in L$, $x^R \in R$, ${x^L}'\in L'$, ${x^R}' \in R'$,
\[\begin{array}{r@{\,}c@{\,}l}
x^L    &< \cut{L'}{R'} &< x^R \\
{x^L}' &< \cut{L}{R}   &< {x^R}'.
\end{array}\]
\end{itemize}
\end{thm}
\begin{definizione}[Canonical representation]
Given $x\in\No$, let 
\[
\begin{aligned}
L^x & := \set{y\in\No: y<x \et y\simple x}\\
R^x & := \set{y\in\No: y>x \et y\simple x}.
\end{aligned}\]
Then, 
\[
x= \cut{L^x}{R^x}.
\]
$\cut{L^x}{R^x}$ is called the \intro{canonical representation} of $x$.
\end{definizione}
\begin{thm}[Inverse cofinality theorem]
Let $x,z\in\No$, $z\simple x$.
Let $x=\cut{L}{R}$ be any representation of $x$.
Then:
\begin{itemize}
\item If $z<x$, there exists $y\in L$ such that $z \leq y <x$.
\item If $z>x$, there exists $y\in R$ such that $z \geq y >x$.
\end{itemize}
\end{thm}
\section{Recursive definitions}
\subsection{Functions of one variable}
\begin{definizione}[Recursive functions]\label{def:recdef}
Let $\ff:\No\to\No$ be a function, $L,R$ be two sets of functions.
I write
\[
\ff = \cut{L}{R}
\]
iff for all $x\in \No$
\begin{equation}\label{eq:recdef}
\ff(x) = \bigcut {\ffL\Paren{x,x^L,x^R,\ff(x^L),\ff(x^R)}}
              {\ffR\paren{x,x^L,x^R,\ff(x^L),\ff(x^R)}},
\end{equation}
where $x^L,x^R$ vary in $L^x$ and $R^x$ respectively, 
and $\ffL,\ffR$ vary in $L$ and $R$ respectively.
\end{definizione}

The formula \eqref{eq:recdef} gives a recursive definition of $\ff$;
in fact, if $\ff(\xo)$ has already been defined for every $\xo$ canonical option of $x$,
it defines uniquely $\ff(x)$ as the simplest element in the cut
\[
\bigcuts {\ffL\Paren{x,x^L,x^R,\ff(x^L),\ff(x^R)}}
      {\ffR\Paren{x,x^L,x^R,\ff(x^L),\ff(x^R)}}_{x^L \in L^x, x^R \in R^x},
\]
if it is non-empty.

The elements $\ffo$ of $L\cup R$ are called \emph{options} of $\ff$;
they are functions with codomain $\No$ and domain classes containing $A\times B$, with
\[\begin{aligned}
A &= \set{(x,y,z)\in\No^3: y<x<z}\\
B &= \set{(\ff(y),\ff(z)): y<z\in\No}.
\end{aligned}\]
I will often use the notations
\[
\ff= \cut{\ffL}{\ffR}
\]
instead of $\ff=\cut{L}{R}$, and
\[
\ffL(x, x^L, x^R)
\]
instead of $\ffL\Paren{x,x^L,x^R,\ff(x^L),\ff(x^R)}$, and similarly for $\ffR$.
\begin{definizione}[Uniform definitions]
The recursive definition $\ff=\cut{\ffL}{\ffR}$ is \intro{uniform} iff the value of
$\ff(x)$ does not depend on the chosen representation of $x$.
This means that: 
\begin{itemize}
\item $\forall x,y,z\in\No$ such that $y<x<z$
\[
\ffL\Paren{x,y,z,\ff(y),\ff(z)}<\ff(x)<\ffR\Paren{x,y,z,\ff(y),\ff(z)}.
\]
\item If $x\in\No$ and $x=\cut{x^L}{x^R}$ is \emph{any} representation of $x$, then
\[
\ff(x) = \cut{\ffL(x,x^L,x^R)}{\ffR(x,x^L,x^R)}.
\]
\end{itemize}
\end{definizione}
If $\Afam$ is a family of functions, and $\ff:\No\to\No$, then $f$ is (\emph{uniformly}) \emph{recursive} over $\Afam$ iff there exist two subsets $L,R$  of $\Afam$ such that $\ff=\cut{L}{R}$ is a (uniform) recursive definition of $\ff$.

Analogous definitions can be given for $\ff$ a function having as domain a convex subclass of $\No$ of the form $\cuts{L'}{R'}$, with $L' < R'$ sub\emph{set}s of $\No$.

If $a=\cut{L}{R}$ then the variable $\ao$, an option of $a$,  will range in $L \cup R$, $a^L$ will be an element of $L$, $a^R$ of $R$.
If do not specify otherwise, $\cut{L}{R}$ will be the canonical representation of $a$, unless I am defining $a$, \ie I am constructing $L$ and $R$.

Similarly, if $f=\cut{L}{R}$, $\ffo$ will range in $L\cup R$, $\ffL$ in $L$, $\ffR$ in $R$.
\subsection{Functions of many variables}\label{sec:manyVars}
In this subsection, $n>0$ is a fixed natural number.
$K$ is the set 
\[
K:=\set{+,-,0}^n\setminus\{(0,\dotsc,0)\}\] 
and $k:=\abs{K}=3^n-1$.
\begin{definizione}
The partial order $\simpleq$ on $\No$ induces a well-founded partial order on $\No^n$ given by
\[
(x_1,\dotsc,x_n) \simpleq (y_1\,\dotsc,y_n) \qquad \text{iff} \qquad
\forall i = 1, \dotsc, n\ x_i\simpleq y_i.
\]
Given $\sigma\in K$ and $\xv,\yv\in\No^n$,
$\yv$ is a $\sigma$-option of $\xv$ iff $\yv\simple \xv$ and for all $i=1,\dotsc,n$
\[\begin{aligned}
y_i & < x_i &\text{iff} && \sigma_i &= -\\
y_i & > x_i &\text{iff} && \sigma_i &= +\\
y_i & = x_i &\text{iff} && \sigma_i &= 0.
\end{aligned}
\]
\end{definizione}
\begin{definizione}
Let $\ff:\No^n\to\No$, $L,R$ be sets of functions.
I write \mbox{$\ff = \cut{L}{R}$} iff for all $\xv\in\No^n$
\[\begin{aligned}
f(x) &= 
  \Splitcut{\ffL\Paren{\xv,\xv^{\sigma(1)},\dotsc,\xv^{\sigma(k)},\ff(\xv^{\sigma(1)}),\dotsc,\ff(\xv^{\sigma(k)})}}
{\ffR\Paren{\xv,\xv^{\sigma(1)},\dotsc,\xv^{\sigma(k)},\ff(\xv^{\sigma(1)}),\dotsc,\ff(\xv^{\sigma(k)})}},
\end{aligned}\]
where $\ffL,\ffR$ vary in $L$ and $R$ respectively, $\sigma$ is a fixed enumeration of $K$, and the $\xv^{\sigma(i)}$ vary among the
$\sigma(i)$-options of $\xv$.
\end{definizione}
For shorthand, I will write $\ffL(x,\xo,\ff(\xo))$
instead of
\[
\ffL\Paren{\xv,\xv^{\sigma(1)},\dotsc,\xv^{\sigma(k)},
     \ff(\xv^{\sigma(1)}),\dotsc,\ff(\xv^{\sigma(k)})}.
\]
Again, the definition of $\ff$ is uniform iff $\ff(\xv)$ does not depend on the chosen representations of the components of $\xv$.
\subsection{Examples}
\begin{example}
Let $\ff(x) := x+1 = \cut{x^L+1,x}{x^R+1}$.
Then, 
\[\begin{aligned}
\ff&=\cut{\ffL_1,\ffL_2}{\ffR}&
\text{where}\\
\ffL_1\Paren{\ff(x^L)} &:= \ff(x^L) & \ffL_2(x)& := x\\
\ffR\Paren{\ff(x^R)} & := \ff(x^R).
\end{aligned}\]
\end{example}
As we can see in the previous example, it might be that some option of $\ff$ does not depend on some of the variables; for instance, $\ffL_1$ does not depend on $x^R$ nor on $\ff(x^R)$.
In this case, I have to impose that $\ff(x)>\ffL_1\paren{\ff(x^L)}$, even for those $x$ that do not have canonical right options.
\begin{example}\label{ex:bsum}
In general, the recursive definition of $\ff(x,y):=x+y$ is
\[
\ff(x,y) = \cut{\ff(x^L,y),\ff(x,y^L)}{\ff(x^R,y),\ff(x,y^R)}.
\]
\end{example}
\begin{example}\label{ex:prod}
Let 
\[\begin{aligned}
\ff(x,y) := xy &= \Splitcut{x^L y + x y^L - x^L y^L, x^R y + x y^R - x^R y^R}
                           {x^L y + x y^R - x^L y^R, x^R y + x y^L - x^R y^L}.
\end{aligned}\]
Then, among the left options of $\ff$ there is 
\[
\ffL_1\Paren{\ff(x^L,y),\ff(x,y^L),\ff(x^L,y^L)} := \ff(x^L,y) + \ff(x,y^L) - \ff(x^L,y^L).
\]
\end{example}
\begin{example}\label{ex:sum}
Let $\ff_0(x),\dotsc,\ff_n(x)$ be recursive functions, $\ff:= \ff_0 + \dotsb + \ff_n$.
Then, 
\[\begin{aligned}
\ffL(x) & = \ff(x) + \ffL_i(x) - \ff_i(x)\\
\ffR(x) & = \ff(x) + \ffR_i(x) - \ff_i(x),
\end{aligned}\]
where $0\leq i \leq n$ and $\ffL_i,\ffR_i$ are left and right options of $\ff_i$.
\end{example}
\begin{proof}
Induction on $n$, using the recursive definition of $x+y$.
\end{proof}
\begin{example}\label{ex:monomial}
Let $a\in\No$, $m \in \Nat$,
$\ff(x):= a x^m$.
Then, 
\[
\ffo(x, x^L, x^R) = ax^m - (a-\ao)(x-x^L)^{\alpha}(x-x^R)^{m-\alpha},
\]
where $0\leq \alpha\leq m$.
Moreover, $\ffo$ is a left option iff $m-\alpha$ is even and $\ao<a$, or $m-\alpha$ is odd and $\ao>a$.
\end{example}
\begin{proof}
Let $x,y\in\No$. The recursive definition of $xy$ implies that an option of $xy$ is of the form
\[
\Opt{xy} = xy - (x-\xo)(y-\yo),
\]
\ie
\[
xy- \Opt{xy} = (x-\xo)(y-\yo).
\]
Therefore, by induction on $m$,
\[
x_1 \dotsm x_n - \Opt{x_1 \dotsm x_m} = (x-\xo_1) \dotsm (x-\xo_m).
\]
In particular, 
\[
\Opt{ax^m} = ax^m - (a-\ao)\Paren{x-(\xo)_1}\dotsm \Paren{x-(\xo)_m},
\]
where $(\xo)_1\dotsm (\xo)_m$ are options of $x$.
Now, I can choose among $(\xo)_1\dotsc,(\xo)_m$ the `best' (\ie the greatest) left option $x^L$, and the `best' (\ie the smallest) right option $x^R$, proving the result.
\end{proof}
\begin{lemma}\label{lem:polynomial}
Let 
\[
\pf(x)=\sum_{0\leq i \leq n} a_i x^i \in \No[x].
\]
Then,
\[
\pfo(x,x^L,x^R) = p(x) - (a_m - \ao_m)(x-x^L)^\alpha(x-x^R)^{m-\alpha},
\]
where $0\leq \alpha \leq m\leq n$.
Moreover, $\pfo$ is a left option iff $m-\alpha$ is even and $\ao<a$, or $m-\alpha$ is odd and $\ao>a$.
\end{lemma}
\begin{proof}
Put together \ref{ex:sum} and \ref{ex:monomial}.
\end{proof}
The relation $\simpleq$ induces a partial order on $\No[x]$.
\begin{definizione}\label{def:polSimple}
Let $\pf,\qf\in\No[x]$.
The polynomial $\qf$ is strictly simpler than $\pf$, in symbols $\qf\simple \pf$, iff
\[\begin{aligned}
\pf &= \sum_{0\leq i\leq n}a_i x^i\\
\qf &= \sum_{0\leq i\leq n}b_i x^i,
\end{aligned}
\]
and there exists a (unique) $m\leq n$ such that
\[\begin{aligned}
b_i & =a_i &i=m+1,\dotsc,n\\
b_m & \simple a_m.
\end{aligned}\]
\end{definizione}
\begin{remark}\label{rem:polOpt}
$\simpleq$ is a well-founded partial order on $\No$, therefore I can do induction on it.
Moreover, if $\pf$ is a polynomial and $\pfo$ is one of its canonical options, as defined in Lemma~\ref{lem:polynomial}, then, as a function of $x$, $\pfo\simple \pf$.
This means that by substituting any value $b ,b'$ for $x^L$ and $x^R$ in $\pfo$ we obtain a polynomial $\qf(x):=\pf(x,b,b')\in\No[x]$ which is strictly simpler than  $\pf(x)$, in the ordering of $\No[x]$.
\end{remark}
\begin{proof}
\[
\pfo(x,x^L,x^R) = \pf(x) - (a_m-\ao_m)(x-x^L)^\alpha(x-x^R)^{m-\alpha}.
\]
As a polynomial in $x$, the $i$-coefficients of $\pfo$ are equal to $a_i$ for $i>m$, 
while the $m$-coefficient is $\ao_m$, which is strictly simpler than $a_m$.
\end{proof}
%

%
\section{Main results}\label{sec:results}
$\DNo$, the Dedekind completion of $\No$, is not a class.
Nevertheless, I can use it as an abbreviation of a well formed formula of \NBG (with a free class variable).

$\BNo$ is the class $\No\cup\set{\pm\infty}$.
In the following definitions, $\Kfield$ is a class, $\leq$ is a linear ordering on it and $\DKfield$ is the Dedekind completion of $(\Kfield, \leq)$.

For the rest of this section, $\Afam$ will be some family of functions on  $\No$.
A function $\ff:\No\to\No$ is really a class, therefore a family of functions is not a class, but only an abbreviation for a well formed formula of \NBG.
\begin{definizione}[Tame]\label{def:tame}
Let $n\geq 0 \in\Nat$, $\ff:\Kfield^{n+1}\to\Kfield$.
$f$ is \emph{tame} iff for every $\dv\in\Kfield^n$ either $\ff(x,\dv)$ is constant, or for every $\zeta\in\DKfield\setminus\Kfield$, $c\in\Kfield$ there exist $a,b\in\Kfield$ such that $a<\zeta<b$ and either
\[\begin{aligned}
\forall x\in (a,\zeta)\ & \ff(x,\dv) > c & \text{or}\\
\forall x\in (a,\zeta)\ & \ff(x,\dv) < c,
\end{aligned}\]
and similarly for $(\zeta,b)$.
\end{definizione}
\begin{definizione}[$\sup$ property]\label{def:sup}
Let $n\geq 0\in\Nat$, $\ff:\Kfield^{n+1}\to \Kfield$.
$\ff$ satisfies the \intro{$\sup$ property} iff $\forall \dv\in\Kfield^n$, $\forall a<b\in\Kfield$,
$\forall c\in \Kfield$, the infima and suprema of the following classes
\[\begin{aligned}
\set{x\in\Kfield: a<x<b \et \ff(x,\dv)\leq c}\\
\set{x\in\Kfield: a<x<b \et \ff(x,\dv)\geq c}
\end{aligned}\]
are in $\Kfield\cup\set{\pm\infty}$.
\end{definizione}
\begin{definizione}[Intermediate value]\label{def:ivp}
A function $\ff:\Kfield^{n+1}\to\Kfield$ satisfies the \intro{intermediate value property} (I.V.P.) iff
$\forall \ev\in\Kfield^n$ $\forall a<b\in \Kfield$ $\forall d\in\Kfield$
\[
\ff(a,\ev)<d<\ff(b,\ev) \to 
\exists c\in \Kfield\ a<c<b\et\ff(c,\ev) = d.
\]
\end{definizione}
Note that the tameness of a given $\ff$ is not a well formed formula of \NBG, because it involves a quantification over elements of $\DNo$, \ie over classes.
Therefore in general it is not possible to speak about the collection of all tame functions in a given collection $\Afam$.
Moreover the theorems involving the tameness of $\Afam$ are actually meta-theorems.

On the other hand, the intermediate value and the $\sup$ properties correspond to well formed formulae of \NBG, because they involve quantifications only on elements of $\No$.
\begin{remark}\label{rem:ivp}
If $\ff:\Kfield\to\Kfield$ satisfies the $\sup$ property and is continuous, then it satisfies the I.V.P..
\end{remark}
\begin{proof}
$c= \sup\set{x\in (a,b): \ff(x)\leq d}$.
\end{proof}
Given a family of functions $\Afam$ and a property $P$ of functions, I say that $\Afam$ satisfies $P$ iff every function in $\Afam$ satisfies $P$; for instance, I could say that $\Afam$ is tame.
\begin{thm}\label{thm:sup}
Suppose that $\Afam$ is tame and satisfies the $\sup$ property.
Let ${\ff:\No\to\No}$ be uniformly recursive over $\Afam$ and tame.
Then, $\ff$ satisfies the $\sup$ property.
\end{thm}
\begin{proof}
Let $a,b,d\in\No$, $a<b$.
Define
\[
\zeta := \sup\set{x\in\No: a<x<b\et \ff(x)\leq d}\in\DNo\cup\set{\pm\infty}.
\]
We have to prove that $\zeta\in\BNo$.
I will prove it by induction on $d$.
Suppose not. Then, $a < \zeta < b$, and, by tameness, \wloG we can suppose that $\ff(x) < d$ in the interval $(a,\zeta)$, while $\ff(x)>d$ in $(\zeta,b)$.
We will construct a $c\in (a,b)$ such that $\ff(c)=d$.

I will give the options of $c$.
First, $a<c<b$, therefore $a$ is a left option, $b$ a right one.
Assume that $\ff=\cut{\ffL}{\ffR}$. Then, $\ff(c)=d$ is equivalent to
\begin{subequations}
  \begin{align}
    \label{eq:sup1}
    d^L&<\ff(c)<d^R\\
    \label{eq:sup2}
    \ffL(c,c^L,c^R)&<d<\ffR(c,c^L,c^R)
  \end{align}
\end{subequations}

\eqref{eq:sup1}
By inductive hypothesis,
\[
c^L := \sup\set{x\in(a,b): \ff(x) \leq d^L}\in\BNo.
\]
Moreover, $c^L \leq \zeta$, and $c^L \neq \zeta$ because $\zeta \notin \BNo$.
Therefore, we can add $c^L$ to the left options of $c$.
Similarly, we find $c^R$ using $d^R$.

\eqref{eq:sup2}
Assume that we have already found some options $c^L, c^R$ of $c$.
Let $U := (c^L,c^R)$.
I will find some `new' options $\cLp,\cRp$ such that $c^L \leq \cLp < \zeta < \cRp \leq c^R$ and if $x\in (\cLp,\cRp)$, then
\[
\ffL(x,c^L,c^R)<d<\ffR(x,c^L,c^R).
\]
Then, I add $\cLp$ and $\cRp$ to the options of $c$, and repeat the process.

Let
\[
\cRp := \inf\set{x\in U: \ffL(x,c^L,c^R)\geq d}.
\]
The $\sup$ property implies that $\cRp\in\BNo$.
If $\cRp= +\infty$, then $\ffL(x,c^L,c^R)\leq d$ in all $U$, therefore I do  not need to add any option to $c$.
Otherwise, $\zeta \leq \cRp \leq c^R$, and $\zeta\neq\cRp$, therefore I can add $\cRp$ to the right options of $c$.

Similarly, let
\[
\cLp := \sup\set{x\in I: \ffR(x,c^L,c^R)\leq d}.
\]
Therefore, at the end of the process we find a $c\in(a,b)$ such that $\ff(c)=d$, a contradiction.
\end{proof}
Note that in the previous theorem I am assuming that $\ff$ is a function of only one variable.
\begin{example}
Conway proves that, with the already defined $<$, $+$ and $\cdot$, $\No$ is a linearly ordered field.
Using the previous theorem, I will show that it is actually real closed.
Conway proves the same thing, but with a quite different technique.
\end{example}
\begin{proof}
An linearly ordered field $\Kfield$ is real closed  iff every polynomial in $\Kfield[x]$ satisfies the I.V.P..
Therefore, by Remark~\ref{rem:ivp}, it is enough to prove that every $\pf(x)\in\No[x]$ satisfies the $\sup$ property.

Moreover, since $\deg p' < \deg p$, the derivative $p'$ is simpler than $p$.
Therefore, by inductive hypothesis, every root of $p'$ is in $\No$, and hence $p$ is tame.

By Remark~\ref{rem:polOpt} every option $\pfo(x,\xo)$ of $\pf(x)$ is simpler than $\pf(x)$, therefore, by induction on $\pf$, it satisfies the sup property.
The conclusion follows from Theorem~\ref{thm:sup}.
\end{proof}
Note that from the proof of Theorem~\ref{thm:sup} it is possible to extract an algorithm to compute the root of a polynomial in $\No[x]$.
This algorithm gives Conway's formula to compute $\unosu{a}$ if the polynomial in question is $ax-1$, and  C.~Bach's formula for $\sqrt{a}$ if we use the polynomial $x^2-a$ instead.
Higher degree polynomials yield quite complicate formulae.
\subsection{Initial substructures of \No}
\begin{definizione}[Initial]\label{def:initial}
Let $S$ be a subclass of $\No$.
$S$ is \intro{initial} iff
\[
\forall x\in S\ \forall y\in\No\ y\simpleq x\to y\in\No.
\]
\end{definizione}
\begin{lemma}\label{lem:initial}
Let $S\subseteq\No$ be an initial subclass of $\No$, $L<R$ be subclasses of $S$,
$x:=\cut{L}{R}$ (if it exists).
Let $z\simple x$.
Then, $z\in S$.
\end{lemma}
Note that I cannot conclude that $x\in S$.
\begin{proof}
\Wlog, I can suppose $z<x$.
Then, by the inverse cofinality theorem, there exists $y\in L$ such that $z\leq y<x$.
So, $z\simpleq y$.

However, $y\in L\subseteq S$ and $S$ is initial, therefore $z\in S$.
\end{proof}
\begin{remark}
The union of an arbitrary family of initial subclasses of $\No$ is initial.
Therefore, given $S\subseteq\No$, I can speak of the maximal initial subclass of $S$;
it is the union of all initial sub\emph{set}s of $S$, therefore it really exists.
\end{remark}
\begin{definizione}[Closure]\label{def:close}
Let $S\subseteq\No$, $\ff:\No^n\to\No$.
$S$ is \intro{closed} under $\ff$ iff
\[
\ff(S^n)\subseteq S.
\]
$S$ is closed  under $\Afam$ iff it is closed under every $\ff$ in $\Afam$.

The \intro{closure} of $S$ under $\Afam$, $\closA{S}$, is the smallest $T\subseteq \No$ closed under $\Afam$ and containing $S$.
\end{definizione}
\begin{thm}\label{thm:initial}
Suppose that for every $T\subseteq\No$ initial, $\closA{T}$ is also initial.
Let $\ff:\No^n\to\No$ be recursive over $\Afam$.
Let $S_1,\dotsc,S_n$ be initial subclasses of $\No$, 
$S:=S_1\times\dotsb\times S_n$.
Then, $\closA{\Paren{S_1\cup\dotsb\cup S_n\cup\ff(S)}}$ is also initial.

Moreover, if every $\ffo$ option of $\ff$ is a function of $\ff(\xo)$ only, 
then $\closA{\paren{\ff(S)}}$ is initial.
\end{thm}
Note that I am not assuming that the definition of $\ff$ is uniform.
\begin{xproof}
Let $U := \closA{\Paren{S_1\cup\dotsb\cup S_n\cup\ff(S)}}$, and
let $T$ be the maximal initial subclass of $U$.
We need to prove that $T = U$.
\begin{claim}\label{cl:closure}
For every $\av\in S$, $b:=\ff(\av)\in T$.
\end{claim}
The proof is by induction on $\av$.
An option of $b$ is
\[
\opt b = \ffo\Paren{\av,\avo,\ff(\avo)}.
\]
Every $\avo$ is strictly simpler than $\av$, therefore $\avo\in S$.
So, by inductive hypothesis, $\ff(\avo)\in T$.
But $S\subseteq T^n$ and, by hypothesis on $\Afam$, $T$ is closed under $\Afam$,
therefore $\opt b\in T$.
Thus, by Lemma~\ref{lem:initial}, $b\in T$.
\begin{claim}
$T=U$.
\end{claim}
It suffices to prove that $T$ contains $S_1, \dotsc, S_n$ and $\ff(S)$, and that $T$ is closed under $\Afam$:
\begin{itemize}
\item $T=\closA{T}$ by hypothesis on $\Afam$.
\item $S_i\subseteq T$, $i=1,\dotsc,n$ because the $S_i$ are initial.
\item $\ff(S)\subseteq T$ by Claim~\ref{cl:closure}.
\end{itemize}
To prove the second point, define $U:=\closA{\paren{\ff(S)}}$, $T$ its maximal initial subclass.
As before, it is enough to prove the following claim:
\begin{claim}
For every $\av\in S$, $b:=\ff(\av)\in T$.
\end{claim}
The proof is by induction on $\av$.
An option of $b$ is
\[
\opt b = \ffo\Paren{\ff(\avo)}.
\]
$\avo\simple \av$, therefore $\avo\in S$.
Thus, by inductive hypothesis, $\ff(\avo)\in T$.
By hypothesis on $\Afam$, $T$ is closed under $\Afam$,
implying that $\opt b = \ffo\Paren{\ff(\avo)} \in T$.
Hence, by  Lemma~\ref{lem:initial}, $b\in T$.
\end{xproof}
\begin{thm}\label{thm:initialCl}
Suppose that for every $T\subseteq\No$ initial, $\closA{T}$ is also initial.
Let $\ff:\No^n\to\No$ be recursive over $\Afam$.
Let $S$ be an initial subclass of $\No$.
Then, $\closAf{S}$ is also initial.
\end{thm}
\begin{proof}
Let 
\[\begin{aligned}
T_0 & := S \\
T_{i+1} & := \closA{\Paren{T_i\cup\ff({T_i}^n)}} \\
T &:=\bigcup_{i\in \Nat} T_i = \closAf{S}.
\end{aligned}\]
Then, by Theorem~\ref{thm:initial} and induction on $i$, $T_i$ is initial for every $i\in\Nat$, therefore $T$ is also initial.
\end{proof}
A different way of presenting the same reasoning is the following.
\begin{proof}
Let $T$ the maximal initial subclass of $\closAf{S}$.
By hypothesis on $\Afam$, $T$ is closed under $\Afam$.
I have to prove that $T$ is also closed under $\ff$.

Let $\av\in T^n$.
I will prove that $\ff(\av)\in T$ by induction on $\av$.
By Lemma~\ref{lem:initial}, it is enough to find $L<R\subseteq T$ such that 
\[w:=\ff(\av) = \cut{L}{R}.\]
An option of $w$ is of the form
\[
\ffo\Paren{\av,\avo,\ff(\avo)},
\]
where $\avo\in\Afam$ and $\avo\simple\av$.
Therefore, $\avo\in T^n$, and, by inductive hypothesis, $\ff(\avo)\in T^n$.
Thus, $\ffo\paren{\av,\avo,\ff(\avo)}\in T$.
\end{proof}
\begin{corollary}\label{cor:initial}
Let $S,U$ be initial subclasses of $\No$. 
Then, the following subclasses of $\No$ are also initial:
\begin{enumerate}
\item The (additive) subgroup generated by $S$.
\item The subring generated by $S$ .
\item $-S :=\set{-x: x\in S}$.
\item $S+U := \set{x+y: x\in S, y\in U}$.
\item The subgroup generated by $SU:=\set{xy: x\in S, y\in U}$.
\end{enumerate}
\end{corollary}
\begin{proof}
For the first two points, apply Theorem~\ref{thm:initialCl}.
The third point is obvious.
For the other two points, apply Theorem~\ref{thm:initial}.
\end{proof}
\begin{example}
It is not true in general that if $S,U$ are initial subgroups of $\No$, then $SU$ is also initial.
For instance, take $S=U$ to be the subgroup generated by $\Zed$ and $\omega$.
Then, $\omega^2 +\omega  = \omega (\omega + 1)\in SU$, but $\omega^2 + 1 \notin SU$.
\end{example}
\begin{corollary}\label{cor:Kx}
Let $\Kfield$ be an initial subring of $\No$.
Let $L<R$ be subsets of $\Kfield$, and let $c:=\cut{L}{R}$.
Then, $\Kfield[c]$ is also an initial subring of $\No$.
\end{corollary}
\begin{proof}
By Remma~\ref{lem:initial}, $\Kfield\cup\set{c}$ is initial, therefore, by Corollary~\ref{cor:initial}, the ring generated by it is also initial.
\end{proof}
\begin{definizione}[Strongly tame]\label{def:axiom5}
A function $\ff: \No^{n+1} \to \No$ is \intro{strongly tame} iff
for all $a<b \in \No$, $\ev \in \No^n$, $d \in \No$ either $\ff(x, \ev)$ is constant, or there exist $\zeta_0, \dotsc, \zeta_{m}\in\DNo$ such that $a = \zeta_0< \dotsb < \zeta_{m} = b$ and for $i = 0, \dotsc, m-1$
\[\begin{aligned}
\forall x\in (\zeta_i,\zeta_{i+1})&\,\ \ff(x,\ev)>d & \text{or}\\
\forall x\in (\zeta_i,\zeta_{i+1})&\,\ \ff(x,\ev)<d.
\end{aligned}\]
\end{definizione}
\begin{definizione}\label{def:closSol}
Let $\ff:\No^{n+1}\to\No$ be strongly tame, $S$ be a subclass of $\No$.
$S$ is \intro{closed under solutions} of $\ff$ iff for all $\ev\in S^n$ either
$\ff(x,\ev)$ is constant, or for every $d\in S$
\[
\forall c\in\No\ \ff(c,\ev)=d\ \Rightarrow\ c\in S.
\]
The \intro{closure} of $S$ under solutions of $\ff$ is the smallest class containing $S$ and closed under $\ff$ and under its solutions.
\end{definizione}
\begin{thm}\label{thm:closSol}
Suppose that $\Afam$ is strongly tame and satisfies the I.V.P..
Let $\ff:\No\to\No$ uniformly recursive over $\Afam$, strongly tame and satisfying the I.V.P..

Suppose that for every $S$ initial subclass of $\No$, the closure of $S$ under solutions of $\Afam$ and under $\ff$ is also initial.
Then, the closure of $S$ under solutions of $\Afamf$ is initial.
\end{thm}
The proof is quite similar to the one of Theorem~\ref{thm:sup}.
\begin{proof}
Let $T$ be the maximal initial subclass of the closure of $S$ under solutions of $\Afamf$.
By hypothesis, $T$ is closed under solutions of $\Afam$ and under $\ff$.
Therefore, it is enough to prove that $T$ is closed under solutions of $\ff$.
If $\ff$ is constant, there is nothing to prove.
Otherwise, let $d\in T$, $c\in\No$ such that $\ff(c) =d$.
I will prove that $c\in T$ by induction on $d$.

I will give options of $c$ in $T$.


Assume that $\ff=\cut{\ffL}{\ffR}$. Then, $\ff(c)=d$ is equivalent to
\begin{subequations}
  \begin{align}
    \label{eq:zero1}
    d^L&<\ff(c)<d^R\\
    \label{eq:zero2}
    \ffL(c,c^L,c^R)&<d<\ffR(c,c^L,c^R)
  \end{align}
\end{subequations}

\eqref{eq:zero1}
Let $\Do$ be an option of $d$.
\[\begin{aligned}
c^L &:= \sup\set{x \in \No: x \leq c \et \ff(x)=\Do}\\
c^R &:= \inf\set{x \in \No: x \geq c \et \ff(x)=\Do}.
\end{aligned}\]
The function $\ff$ is  strongly tame and satisfies the I.V.P., thus $c^L,c^R\in\BNo$ and $c^L \leq c \leq c^R$.
By the I.V.P., $c^L < c < c^R$.
Moreover, $c^L, c^R\in T \cup {\pm\infty}$ by inductive hypothesis.
By the I.V.P., $\ff(x) - \Do$ does not change sign in $(c^L,c^R)$; 
in particular, if $\Do = d^L<d$, $\ff(x)>d^L$ in $(c^L,c^R)$, and similarly for $\Do=d^R$.
Hence, I can take $c^L$ and $c^R$ as left and right options of $c$.

\eqref{eq:zero2}
Suppose that I have already found $c^L,c^R$ `old' options of $c$.
Let $\ffo$ be an option of $\ff$, say $\ffo=\ffL<\ff$.
I will construct $\cLp,\cRp$ `new' options of $c$ such that for $x\in(\cLp,\cRp)$ $\gf(x):=\ffL(x,c^L,c^R)<d$.
Let
\[\begin{aligned}
\cLp &:= \inf\set{x\in\No: x\leq c \et \gf(x) = d}\\
\cRp &:= \sup\set{x\in\No: x\geq c \et \gf(x) = d}.
\end{aligned}\]
$\gf$ is  strongly tame, therefore $\cLp,\cRp\in\BNo$ and $\cLp<c<\cRp$.
By the I.V.P., $\gf(x)-d$ does not change sign in $(\cLp, \cRp)$; 
consequently, $\gf(x)<d$ in $(\cLp,\cRp)$.
So, I can take $\cLp$ and $\cRp$ as left and right options of $c$.

At the end of the process, I obtain $L<R\subseteq T$ such that $c=\cut{L}{R}$, so, by Lemma~\ref{lem:initial}, $c\in T$.
\end{proof}
\begin{corollary}\label{cor:initReal}
The real closure of an initial subring of $\No$ is initial.
More in general, if $S\subseteq\No$ is initial, then the smallest real closed subfield of $\No$ containing $S$ is initial.
\end{corollary}
\begin{proof}
$\Raz$ is an initial subfield of $\No$, therefore $\Raz\cup S$ is also initial.
Therefore, by Corollary~\ref{cor:initial}, $\Kfield$, the subring generated by it, is also initial.

It remains to prove that the real closure of $\Kfield$ is initial.
Apply the Theorem~\ref{thm:closSol} and induction on $\No[x]$.
\end{proof}
If instead of considering \emph{all} polynomials in $\No[x]$, we consider only the polynomials of degree up to a fixed degree $n$, the previous corollary is still valid, with the same proof.
In particular, taking $n=1$, we can conclude the following:
\begin{corollary}\label{cor:initField}
Let $S$ be an initial subclass of $\No$.
Then, the subfield of $\No$ generated by $S$ is initial.
\end{corollary}
The following theorem was also proved in \cite{ehrlich} with a different method.
\begin{thm}\label{thm:isomRCF}
Let $\Kfield$ be a real closed field and a proper set.
Then, $\Kfield$ is isomorphic to an initial subfield of $\No$.
\end{thm}
\begin{proof}
If $\Kfield\isom\Raz$, it is true.

Assume that $\Ffield$ is a real closed initial subfield of $\No$, and $\Kfield$ is (isomorphic to) the real closure of $\Ffield(a)$ for some $a$ transcendental over $\Ffield$.
Let $\cuts{L}{R}$ be the cut determined by $a$ over $\Ffield$.
For any $c\in \cuts{L}{R}$, $\Ffield(c)$ is isomorphic to $\Ffield(a)$, and its real closure is isomorphic to $\Kfield$.
Moreover if $c=\cut{L}{R}$ then, by Corollary~\ref{cor:Kx}, $\Kfield[c]$ is initial, therefore, by Corollary~\ref{cor:initReal}, its real closure is also initial.

In general, let $(c_\beta)_{\beta<\alpha}$ be a transcendence basis of $\Kfield$ over $\Raz$, for some $\alpha\in\No$.
Let $\Kfield_0$ be the real closure of $\Raz$, and define
$\Kfield_\beta$ to be the real closure of $\Raz[c_i: i<\beta]$ for $\beta<\alpha$, \ie $\Kfield_0:=\Raz$, and, for $0<\beta\leq\alpha$
\[
\Kfield_\gamma := \begin{cases}
\text{the real closure of }\Kfield(c_\gamma) &\text{if } \beta = \gamma+1\\
\underset{\gamma<\beta}{\bigcup}\Kfield_\gamma &\text{if } \beta \text{ is limit.}
\end{cases}\]
In particular, $\Kfield_\alpha := \Kfield$.
By the previous case and induction on $\beta$, each $\Kfield_\beta$ is isomorphic to an initial subfield of $\No$, and the conclusion follows.
\end{proof}
It is not true that every ordered field (which is also a set) is isomorphic to an initial subfield of $\No$.
For instance, take $\Kfield := \Raz(\sqrt{2}+\unosu{\omega})\subset\No$.
Suppose, for contradiction, that there exists an isomorphism of ordered fields  $\psi$ between $\Kfield$ and an initial subfield of $\No$.
Let $z=\psi(\sqrt{2} + \unosu{\omega})$.
Then, $\sqrt{2}\simple z$, but $\sqrt{2}\notin\psi(\Kfield)$.

For more on the subject of initial embeddings of fields, see~\cite{fornasiero2}.
%

%
\section{Examples and the concatenation function}
\begin{example}
Let $\ff(x):=\cut{x-1}{x+1}$.
The image of $\ff$ is the class of omnific integers.
$\ff$ satisfies the $\sup$ property, but not the I.V.P..  
$\ff(x) = 0$ for $x\in(-1,1)$, $\ff(x) = \alpha$ for $x\in[\alpha,\alpha+1)$, $\alpha\in\On$, etc.
\end{example}
\begin{example}
Let 
\[
\ff(x):= [x]-x = \cut{-1, [x]-x^R}{1,[x]-x^L}.
\]
Then, $\ff$ is not tame.
\end{example}
\begin{proof}
Consider the cut $\omega\concat{-\infty}$ between positive finite numbers and infinite numbers.
$\ff(x)$ changes sign infinitely many times in every neighbourhood of this cut.
\end{proof}
\begin{example}
Let 
\[\begin{aligned}
\ff(x) &:= \Cut{-\abs{x}}{\ffR_1(x^L), \ffR_2(x^R)} & \text{with}\\
\ffR_1(z) &:= \begin{cases}
0 & \text{iff } z \geq 0\\
2 & \text{iff } z < 0
\end{cases}\\
\ffR_2(z) &:= \begin{cases}
0 & \text{iff } z \leq 0\\
2 & \text{iff } z > 0.
\end{cases}
\end{aligned}\]
Therefore,
\[
\ff(x) = 
\begin{cases}
\cut{-\abs{x}}{0}<0 & \text{iff } x\neq 0\\
1 & \text{iff } x=0
\end{cases}
\]
\end{example}
\begin{definizione}
Let $x,y>0\in\No$.
\begin{itemize}
\item $x\simeq y$ iff $\frac{x-y}{x}$ is infinitesimal.
\item $x\stless y$ iff $x<y$ and $x\nsimeq y$.
\item $x\perp y$ iff $x\nsimpleq y$ and $y\nsimpleq x$.
\end{itemize}
\end{definizione}
\begin{example}
Let $d:=\dueterzi$.
For $a\in\Real$, let
\[\begin{aligned}
\ffL_a(z) &:= \begin{cases}
z & \text{if } z\leq a\\
a & \text{if } z\geq a
\end{cases}\\
\ffR_a(z) &:= \begin{cases}
z & \text{if } z\geq a\\
a & \text{if } z\leq a.
\end{cases}
\end{aligned}\]
Let 
\[
\ff(x):=\cut{\ffL_r(x^L): d > r\in\Real}{\ffR_s(x^R): d < s\in\Real}.
\]
Then, $\ff$ is uniformly recursive.
Moreover,
\[\begin{aligned}
\ff(x) &= x & \text{if } & x\stless d\\
\ff(x) &= d & \text{if } & x\simeq d\\
\ff(x) &= x & \text{if } & d\stless x.
\end{aligned}\]
In particular, $\ff$ does not satisfy the $\sup$ property, because
\[
\sup\set{x\in\No: 0<x<1\et \ff(x)\leq d} \notin \No
\]
\end{example}
In the previous examples, the definition of $\ff$ is uniform.

A surreal number $z$ can be considered as a function from $\lt(z)$ into $\set{+,-}$, \ie as a sequence of pluses and minuses, the \emph{sign sequence} of $z$.
Therefore, instead of $1$ we can write the corresponding sequence $+$, instead of $\unmezzo$, we can write $+-$, etc.

An element $\zeta$ of $\DNo$ has also a unique sign sequence, of length $\On$ iff $\zeta \notin \No$.
On the other hand, not every sign sequence corresponds to a element of $\DNo$; for instance, $+\infty$, the sequence of pluses of length $\On$, is not in $\DNo$.
\subsection{The concatenation function}
\begin{definizione}[Concatenation]
Let $x,y\in\No$.
The \emph{concatenation} of $x,y$, noted with $x\concat y$, is given by the sign sequence of $x$ followed by the sign sequence of $y$.
\end{definizione}
The recursive definition of $x\concat y$ is
\[
x\concat y = \cut{x^L,x\concat y^L}{x^R, x\concat y^R}.
\]
This definition is \emph{not} uniform; while I can choose \emph{any} representation of $y$, I must take the canonical representation of $x$ (but see also \cite{keddie}).
In the following, I will prove that, given some hypothesis on $\Afam$, it is never possible to find a uniform recursive definition of $x\concat y$  over $\Afam$.

Let 
\[
\ff(x):= x\concat 1 = \cut{x}{x^R}.
\]
By definition, $\ff(x)<0$ for $x<0$, while $\ff(x)>0$ for $x\geq 0$, therefore
$\ff$ does not satisfy the I.V.P..
Moreover, $\ff$ is injective.
\begin{table}[htbp]
  \centering
  \begin{tabular}[b]{|c|ccccccl|}
    \hline
    \raisebox{0ex}[2.5ex][1.3ex]{}
    $x$          & 0 & $\alpha$     & $-\alpha-1$         &
    $-\omega$   & $\unosu{2}$   & $\unosu{\omega}$ &
    \begin{small}$\alpha\in\On$.\end{small}\\
    \hline
    \raisebox{0ex}[2.5ex][1.5ex]{}
    $x\concat 1$ & $1$ & $\alpha+1$ & $-\alpha-\unosu{2}$ &
    $-\omega+1$ & $\smallfrazione{3}{4}$ & $\smallfrazione{2}{\omega}$ & \\
    \hline
  \end{tabular}
  \caption{Some values of $x\concat 1$}
  \label{tab:concat1}
\end{table}

\begin{lemma}
Let $x\in\No$.
Then:
\begin{itemize}
\item $\ff(x) > x$.
\item $x \simple \ff(x)$
\item If $x < y$ and $y \simple x$, then $\ff(x) < y$.
\end{itemize}
\end{lemma}
\begin{proof}
Obvious.  
\end{proof}
\begin{lemma}\label{lem:concatTame}
Let $c,d\in\No$. Then, there exists $a,b \in\BNo$ such that  $a<c<b$ and
\[\begin{aligned}
\forall x\in (a,c) &\ \ff(x)<d & \text{or}\\
\forall x\in (a,c) &\ \ff(x)>d,
\end{aligned}\]
and similarly for $(c,b)$.
\end{lemma}
\begin{proof}
There are three cases, according to $d\lesseqgtr c$.
\begin{itdescriptionD}
\item[$d<c$]
Let $a:=\cut{d}{c}$, $b:=+\infty$.
Then, for every $x>a$ $\ff(x)>x>a>d$.
\item[$d=c$]
Let $a:=c\concat -$, $b:=+\infty$.
For every $x\in (a,c)$ $c\simple x$, therefore $\ff(x)<c$.
For every $x>c$ $\ff(x)>x>c$.
\item[$d>c$]
Let $a:=c\concat -$.
If $a<x<c$, then $c\simple x$, therefore $\ff(x)<c<d$.

Assume that $c\perp d$.
Let $b:=\cut{c}{d}$.
If $c<x<b$, then $b\simple x$, therefore $\ff(x)<b<d$.

Assume that $c\simple d$.
Let $b:=\cut{c}{d}$.
If $c<x<b$, then $b\simple x$, therefore $\ff(x)<b<d$.

Assume that $d\simple c$.
Let $b:=d$.
If $c<x<b$, then $d\simple x$, therefore $\ff(x)<d$.\qedhere
\end{itdescriptionD}
%
%
%
%
%
\end{proof}
\begin{remark}\label{rem:23}
Let $d := \dueterzi = +-+-+-\dotso$,
\begin{align*}
a_0 &:= 0,\\
a_1 &:= +- = \unosu 2,\\
a_2 &:= +-+- = \nicefrac{5}{8},\\
a_3 &:= +-+-+-,\\
&\dotso,\\
a_{n+1} &:= a_n\concat{+-},
\intertext{and} 
b_n &:= a_n\concat + \ =\ +-\dotso+-+\\
c_n &:= a_n\concat - \ =\ +-\dotso+--.
\end{align*}
Then, for every $n\in\Nat$, 
\[
c_n <a_n <c_{n+1} < d < b_n.
\]
Moreover, $\ff(a_n) = b_n > \ff(d) > d$, while $\ff(c_n) < a_n < d$.
Besides,
\[
d =\cut{a_n}{b_n}_{n\in\Nat}
\]
is the canonical representation of $d$.
Moreover for every $x\in\No$
\[
d\simpleq \ff(x) \Leftrightarrow d \simpleq x.
\]
Finally, for every $x\in\No$ $d\simpleq x$ iff $d-x$ is infinitesimal.
\end{remark}
\begin{remark}
$\ff$ is not continuous.
In fact,
\[
\lim_{\substack{x\to 0 \\ x\neq 0}}\ff(x) =0,
\]
while $\ff(0)=1$.
\end{remark}
\begin{lemma}
$\ff$ is not tame.  
\end{lemma}
\begin{proof}
Let $d,a_n,b_n,c_n$ be as in Remark~\ref{rem:23}.

If $a<\zeta$, there exists $n\in\Nat$ such that $a<a_n<c_{n+1}<\zeta$.
Therefore, $\ff(x)-d$ changes sign infinitely many times in every left neighbourhood of $\zeta$.

However, note that in every infinitesimal neighbourhood of $d$, $\ff(x)\geq d$ iff $x\geq d$, and $\ff(x)\geq\ff(d)$ iff $x\geq\ff(d)$ or $x=d$.
\end{proof}
\begin{lemma}
$\ff$ does not satisfy the $\sup$ property.
\end{lemma}
\begin{proof}
Let $d,a_n,b_n,c_n$ be as in Remark~\ref{rem:23}, and let $a:=0$.

Then, $\ff(a_i)>d$, while $\ff(c_i)<d$.
Let
\[
\zeta := \sup\set{x\in\No: a<x<d \et \ff(x)\geq d}.
\]
Therefore, 
\[
\zeta\geq\sup\set{a_i:i\in\Nat}.
\]
On the other hand, if $d\simpleq x$ and $x<d$, then $\ff(x)<d$, therefore
\[
\zeta \leq \inf\set{x<d: d\simpleq x},
\]
so $\zeta = \inf\set{x<d: d\simpleq x}\in\DNo\setminus\No$.
\end{proof}
\begin{lemma}
If $\Afam$ is a family of functions strongly tame and satisfying the I.V.P., then $\ff$ cannot be uniformly recursive over $\Afam$.
\end{lemma}
\begin{proof}
Suppose not, \ie that $\ff=\cut{\ffL}{\ffR}$, with $\ffo\in\Afam$.

Let $d,a_n,b_n,c_n$ be as in Remark~\ref{rem:23}, and let $\zeta$ as in the previous proof.

I will show that there exists $c\in\No$ such that $\ff(c) = d$, which is clearly impossible. I will give the options of $c$.

$\ff(c) = d$ is equivalent to:
\[\begin{aligned}
d &\simpleq \ff(c)\\
\ffL(c,c^L,c^R) &<d < \ffR(c,c^L,c^R).
\end{aligned}\]

By Remark~\ref{rem:23}, $d\simpleq \ff(c)$ is equivalent to $d\simpleq c$, and $d = \cut{a_n}{b_n}_{n\in\Nat}$.
\\
Therefore, it is necessary and sufficient to add $a_n$ to the left options of $c$ and $b_n$  to its right ones for every $n\in\Nat$ to ensure that $d\simpleq \ff(c)$.

Let $c^L,c^R$ be `old' options of $c$ such that $c^L<\zeta<c^R$, $\ffL$ be a left options of $\ff$.
I will find $\cLp,\cRp$ `new' options of $c$ such that $\cLp<\zeta<\cRp$ and
\begin{equation}\label{eq:concatSup}
\forall x\in(\cLp,\cRp)\ \ffL(x,c^L,c^R)<d.
\end{equation}
Let $\gf(x) := \ffL(x,c^L,c^R)$,
\[\begin{aligned}
\cLp &:=\sup\set{x\in\No: x<\zeta \et \gf(x)= d}\\
\cRp &:=\inf\set{x\in\No: x>\zeta \et \gf(x)= d}.
\end{aligned}\]
$\Afam$ is strongly tame, therefore the previous $\sup$ and $\inf$ are actually a $\min$ and a $\max$, unless $\gf(x)$ is constant.

For every left neighbourhood $J$ of $\zeta$ there is $x\in J$ such that $\ff(x)< d$.
Moreover, $\gf(x) < \ff(x)$ in $(c^L,c^R)$.
Therefore, if $\gf(x)$ is constant, then $\gf(x) < d$, so $\cLp = -\infty$ and $\cRp = +\infty$.

Otherwise, $\cLp, \cRp \in \BNo$, while $\zeta \notin \BNo$, so $\cLp < \zeta < \cRp$.
By the I.V.P., the sign of $\gf(x) - d$ is constant in $U := (\cLp,\cRp)$.
Again, in every left neighbourhood of $\zeta$ there is a $x$ such that $\ff(x)<d$, therefore $\gf(x)<d$ in $U$.

Proceed similarly for $\ffR$.

At the end of the process we have constructed a $c\in\No$ such that \mbox{$\ff(c)=d$}, a contradiction.
\end{proof}
%

%
%
\nocite{alling,conway,gonshor,knuth,dries,driesEhrlich,mythesis,
ehrlich,keddie,allingEhrlich:abs,allingEhrlich:alt,kuhlmann}
\bibliography{surreals}
\bibliographystyle{abbrv}
%

%
\end{document}